# COMPUTING SHAPE DNA USING THE CLOSEST POINT METHOD

RACHEL HAN

**Abstract.** We demonstrate an application of the closest point method where the truncated spectrum of the Laplace–Beltrami operator of an object is used to identify the object. The effectiveness of the method is analyzed as well as the default algorithm, 'eigs', in MATLAB which computes the eigenvalues of a given matrix. We also cluster "similar" objects via multi-dimensional scaling algorithm and empirically measure its effectiveness.

**1. Introduction.** The Laplace–Beltrami operator is widely used in geometric modelling and computer graphics for applications such as smoothing, segmentation and registration of 2D or 3D shapes [5], [8]. Another novel application is characterizing shapes by extracting their 'fingerprints' or 'Shape-DNA', as first introduced by Reuter, Wolter and Peinecke [6]. Storing and processing the 'Shape-DNA's enables fast retrieval and identification in a database of shapes and has potential applications in areas such as machine learning. The full spectrum of Laplace–Beltrami on a surface is able to identify distinct shapes because it contains intrinsic information of the Riemann manifold, such as volume and surface area [6]. Also, the identification via the the Shape DNA is robust since it is isometry invariant (isometries include rotation, translation and reflection) and independent of parametrization, reducing the preprocessing of the shapes. In this project, the robustness is improved further via the closest point method [7] which is able to represent surfaces without any parametrization. The Laplace-Beltrami operator on curved surfaces is analogous to the standard Laplacian operator. For example, it models diffusion on a surface. Using the work previously done by Macdonald, Brandman and Ruuth [4] and Von Glehn, März and Macdonald [9], we use the closest point method and MATLAB to discretize and solve the given eigenvalue problem on various surfaces and assess its numerical errors. Furthermore, we analyze the effectiveness of the default MATLAB 'eigs' algorithm and compare its performance with our indirect approach. Finally, we use multidimensional scaling plots to represent the similarities between the Shape DNA's and analyze the results.

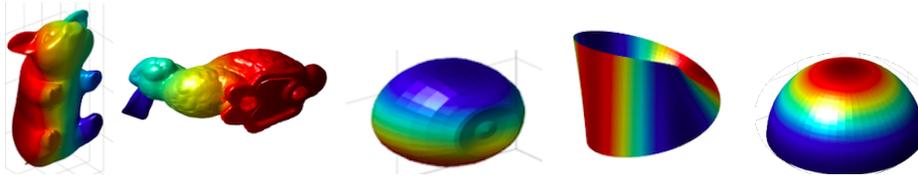

Fig. 1: Eigenfunctions on pig, Stanford bunny, apple, mobius strip and hemisphere surfaces

**2. The Laplace–Beltrami Eigenvalue Problem.** To compute the spectrum of a Laplace-Beltrami operator $\Delta_s$ on a surface $S$, we need to find the eigenpairs $(\lambda, u)$ to the functional problem:

$$-\Delta_s u(\boldsymbol{x}) = \lambda u(\boldsymbol{x}) \text{ for } \boldsymbol{x} \in S,$$

where the *eigenvalue* $\lambda$ is a scalar number and the *eigenfunction* $u$ is a function defined on $S$. This eigenvalue problem follows from solving the wave equation where the vibrating membrane is fixed under a certain boundary condition, via separation of variables:

$$\frac{\partial^2 u}{\partial t^2} = c^2 \nabla^2 u.$$

The eigenfunctions can be interpreted as vibration patterns (standing waves) of the domain, and the eigenvalues are associated with the frequency of the waves. The paper by Kac, *Can One Hear the Shape of a Drum?* [2], poses the question whether one can identify the shape given the eigenvalues of the Laplacian operator. In fact, two different membranes can give an identical set of eigenvalues as shown through isospectral spectral domains constructed by Gordon et al.[1]. However, it is known that membranes with the same area and "connectivity" yield the same set of eigenvalues [2]. Given that we can infer some information of the surfaces using the spectra, if we are given a frequency of a membrane, can we guess which shape of the drum it is most "similar" to? Also, in what ways are the surfaces "similar"? In this paper, we attempt to classify surfaces using the L-B spectra by automatically clustering them, where similarities are optimally determined through a clustering algorithm. We measure the success of the clustering empirically, by observing how close the generated clusters are from the expected clustering based on intuitive characteristics of the surfaces such as roundness of the surface and boundaries. To do this, we require a notion of Shape DNA [6]:

DEFINITION 2.1. *Shape DNA: list of eigenvalues of the L-B operator, scaled by the first non-zero eigenvalue.*

The non-zero scaling factor enables the Shape DNA to be invariant to scaling.

We show a simple analytic calculation of the spectrum on a circle. The eigenvalue problem for the Laplacian in polar coordinates can be formulated as follows:

$$\tag{1} u_{rr} + \frac{1}{r} u_r + \frac{1}{r^2} u_{\theta\theta} + \lambda u = 0$$

Since the Laplacian is defined on a circle of radius 1, $r = 1$ is a constant and $u$ is dependent on $\theta$ only. First we study the boundary value problem of the following ODE with a homogeneous periodic boundary condition:

$$\tag{2} \begin{aligned} u'' + \lambda u &= 0 \\ u(0) = u(2\pi),\ u'(0) &= u'(2\pi) \end{aligned}$$

The solution is a linear combinations of sines and cosines:

$$\tag{3} u(\theta) = A\cos(\sqrt{\lambda}\theta) + B\sin(\sqrt{\lambda}\theta)$$

Applying the Dirichlet boundary condition when $\lambda \geq 0$,

$$\begin{aligned} u(0) = u(2\pi) &= 0 \\ \Rightarrow A = 0 \Rightarrow u(\theta) &= B\sin(\sqrt{\lambda}\theta) \\ B\sin(\sqrt{\lambda}2\pi) &= 0 \end{aligned}$$

$$\tag{4} \begin{aligned} &\Rightarrow \lambda \neq 0 \\ \lambda_n &= \left(\frac{n}{2}\right)^2, n = 1, 2, ... \\ u_n &= \sin\left(\frac{n\theta}{2}\right) \end{aligned}$$

When $\lambda = 0$ and $\lambda \leq 0$, it yields a trivial solution so is not an eigenvalue.

With a Neumann boundary condition,

(5)
$$\begin{aligned}
u'' + \lambda u &= 0 \\
u'(0) = u'(2\pi) &= 0 \\
B = -A\sin(\sqrt{\lambda}2\pi) + B\cos(\sqrt{\lambda}2\pi) &= 0 \Rightarrow B = 0 \\
\Rightarrow -A\sin(\sqrt{\lambda}2\pi) &= 0 \\
\lambda_n &= \left(\frac{n}{2}\right)^2, n = 1, 2, ... \\
u_n &= \cos\left(\frac{n\theta}{2}\right)
\end{aligned}$$

Additionally, when $\lambda = 0$, $u(\theta) = C\theta + D$, it yields a non-trivial eigenfunction $u(\theta) = D$, which is a constant function. We let $u_0 = 1$ with the corresponding eigenvalue $\lambda_0 = 0$.

**3. Closest Point Method.** To solve for the eigenvalues for the Laplace-Beltrami operator on a given surface, we use the closest point method [7] to numerically represent the problem. Consider the following surface PDE:

(6)
$$\begin{aligned}
u_t &= \Delta_s u \\
u_0 &: S \mapsto \mathbb{R}
\end{aligned}$$

Given such PDE on the surface, for each Cartesian grid point in a narrow band near the surface, the closest point extension operator copies the function value from the point on the surface that is closest to the grid point. The narrow band near the surface is populated with the function values defined on the surface. The formal definitions are given below.

DEFINITION 3.1. *[7] Closest Point Function: Given a surface $S$, $cp(x)$ refers to a (possibly non-unique) point belonging to $S$ which is closest to $x$.*

DEFINITION 3.2. *[7] Closest Point Extension Operator. Let $S$ be a smooth surface in $\mathbb{R}^d$. The closest point extension of a function $u : S \to \mathbb{R}$ to a neighborhood $\Omega$ of $S$ is the function $v : \Omega \to \mathbb{R}$ defined by $v(x) = u(cp(x))$. $E$ is the operator which maps $u$ to $v$: $v(x) = u(cp(x)) \iff v = Eu$.*

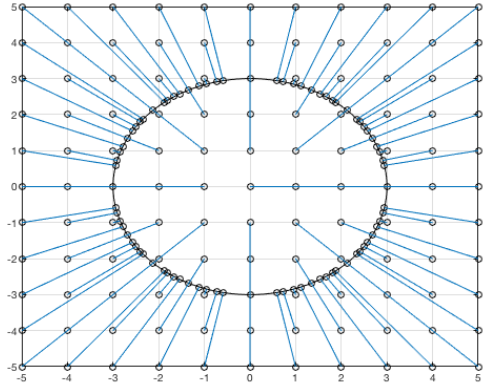

Fig. 2: Closest Point Function on a Circle

We hope to use finite difference schemes to solve the Laplace-Beltrami eigenvalue problem in the Cartesian embedding space. By the gradient principle [7] of the closest point method, the intrinsic surface gradient $\nabla_s u(\boldsymbol{x})$ is same as $\nabla v(\boldsymbol{x})$ because the closest point extension of $u$ is constant in the normal direction. Also, by the divergence principle, $\nabla_s \cdot \boldsymbol{v}(\boldsymbol{x})$ and $\nabla \cdot \boldsymbol{v}(\boldsymbol{x})$ are equal, where $\boldsymbol{v}(\boldsymbol{x})$ is tangent to the surface and all surfaces displaced by a fixed distance from $S$. The principles give us a method of solving a variety of differential equations on surfaces, including the Laplace-Beltrami operator which is the surface divergence of the gradient:

THEOREM 3.3. *[7] Let $S$ be a smooth closed surface in $\mathbb{R}^d$ and $u : S \to \mathbb{R}$ be a smooth function. Assume the closest point function $cp(x)$ is defined in a neighborhood $\Omega \subset \mathbb{R}^d$ of $S$. Then*
$$\Delta_s u(x) = \Delta v(x) \text{ for } x \in S.$$

The problem is now in the Cartesian band surrounding the surface, and its solution to the differential equation must be consistent with the solution of the same equation on the surface ($u = v|_s = Eu|_s$). The new formulation of the problem given in (6) is as follows:

$$
\begin{aligned}
v_t &= E\Delta v \\
v &= Ev \\
v_0 &= Eu_0
\end{aligned}
\tag{7}
$$

The method-of-lines penalty approach [9] is an equivalent formulation that is numerically favourable in practice, with a penalty parameter $\gamma$ which controls the constraint, $v = Ev$.

$$
\begin{aligned}
v_t &= E\Delta v - \gamma(v - Ev) \\
v_0 &= Eu_0
\end{aligned}
\tag{8}
$$

Using the Lagrange interpolating polynomial of appropriate degrees, we can approximate the extension operator $E$ from the neighbouring grid points and discretize the L-B operator, denoted as M:

$$[9] M = E_1 L - \gamma(I - E_3). \tag{9}$$

$E_1$ and $E_3$ are the closest point extension matrices with degree 1 and 3 respectively, chosen experimentally. $L$ is the discretized Laplacian, the standard Cartesian centered approximation with second order of accuracy, using the five point stencil [3].

**4. Implementation.** Once $M$ in (1) is computed, we can use the MATLAB command 'eigs' to compute the eigenvalues and the corresponding eigenvectors ($u$ in vector form). The resulting matrix $M$ is a sparse, nonsymmetric matrix. MATLAB uses Arnoldi iteration to compute the eigenvectors. The idea is similar to power iteration, where we start with a random vector $v$, and iteratively compute $Mv, M^2v, M^3v, ..., M^{n-1}v$. These form a Krylov matrix:

$$\begin{bmatrix} b & Mv & M^2v & \cdots & M^{n-1}v \end{bmatrix}.$$

Arnoldi iteration orthogonalizes this column space through Gram-Schmidt orthogonalization.

This sequence converges to the largest eigenvector $v_{max}$ corresponding to the largest eigenvalue. However, we want the smallest 50 eigenvalues, which forces us to compute the largest eigenvectors for the inverse of M. This requires us to solve the subproblem $Mv_n = v_{n-1}$ every step of the Arnoldi iteration. The default method of solving this linear system in MATLAB is the backslash operator. Since $M$ is a nonsymmetric matrix, the LU solver is used. This direct method does not take advantage of the sparsity of M, and runs out of memory for $\Delta x = 0.0125$ for solving the eigenvalue problem on a sphere of radius 1. To resolve finer grids for accurate computation using smaller memory space, we implemented an iterative solver using GMRES to solve the subproblem of solving $Mv_n = v_{n-1}$. Incomplete L-U factorization was performed for preconditioning, which was essential for the solution to converge in a reasonable amount of time.

The modified iterative solver was tested on 16.43 GB of RAM on a Linux machine with 4x Intel Core Processor (Haswell).

| $\Delta x$ | No Precondition | ILU |
|---|---|---|
| 0.05 | 1047.5 | 132.7 |
| 0.1 | 139.8 | 21.591 |
| 0.2 | 29.8 | 7.5133 |

This approach is slower than the default direct approach. However, we see that there is a trade off between runtime and memory for the finer values of $\Delta x$. The percentage of RAM usage of each solver method was compared. The direct solver uses almost all the available memory when $\Delta x = 0.2/16$, whereas the indirect solver exhibits a slower increase.

**5. Numerical Experiments.**

**5.1. Convergence Studies.** The error in computing the truncated spectrum ($\lambda_n$) of the L-B operator using the closest point method is $C(n)\Delta x^2$; second order convergence as $\Delta x \to 0$. The constant $C(n)$ would grow larger as $n$ increases as smaller $\Delta x$ will be required to resolve the higher eigenvalues as the eigenfunctions become more oscillatory. Errors between the exact and computed eigenvalues of a unit sphere were studied in $\|.\|_\infty$ and $\|.\|_2$. The red dotted lines indicate slope 2 in a log-log scale plot.

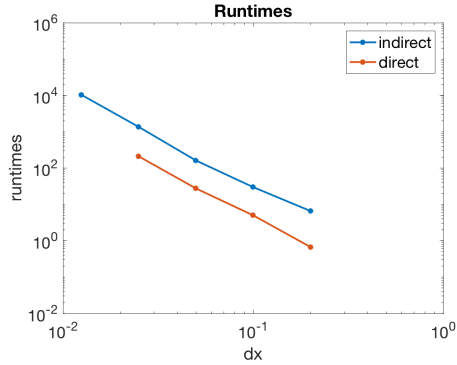
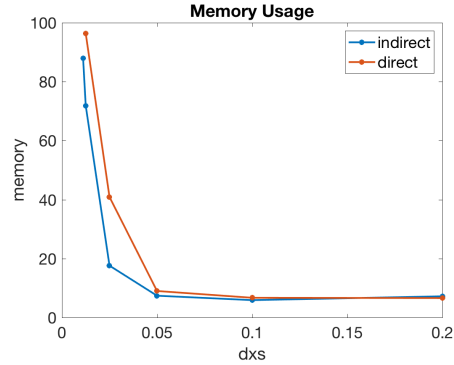

Fig. 3: Runtime

Fig. 4: Memory Usage

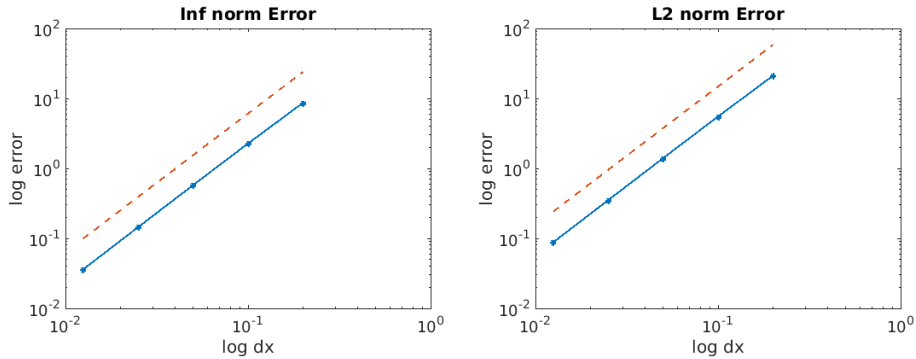

Fig. 5: Error in Computed Spectrum of a Sphere

**5.2. Rotational Invariance.** We also test the rotational invariance, which is one of the isometry properties of the spectrum. We rotate open and closed surfaces which is represented via the closest point method by $\frac{\pi}{5}, \frac{\pi}{4}, \frac{\pi}{3}$ and $\frac{\pi}{2}$ and take the euclidean difference between the non-rotated spectrum and the rotated one. As before, we attempt to verify the second order convergence of errors as $\Delta x \to 0$ because the closest point method of computing the L-B spectrum is second order. Below are the rotational invariance checks for an ellipsoid and a hemisphere.

We observe second order convergence of errors for the ellipsoid, but the errors of the hemisphere show approximately first order convergence. We suspect that this is due to the additional numerical errors associated with the imposed boundary conditions via CPM [7]. Second-order accuracy can be obtained by a minor modification of the extension operator [4].

**5.3. Scaling Invariance.** As mentioned previously, we scale the eigenvalue spectra by the first non-zero vector to ensure the Shape DNA is robust to different scales of shapes studied. For example, as $\Delta x \to 0$, a sphere of radius 1 will have the same Shape DNA as the sphere of radius 2. The error convergence plot over $\Delta x$ was generated using a torus of minor radius 0.5 and major radius 1, and a torus of minor radius 1 and major radius 2.

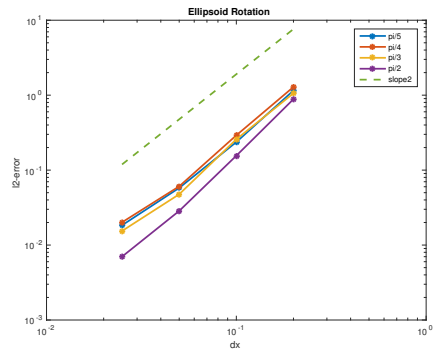

Fig. 6: Ellipsoid

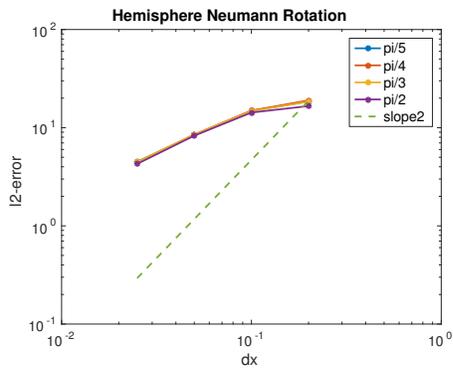

Fig. 7: Hemisphere Neumann BC

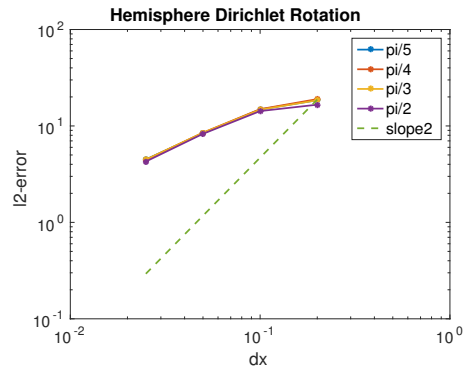

Fig. 8: Hemisphere Dirichlet BC

Rotational Invariance

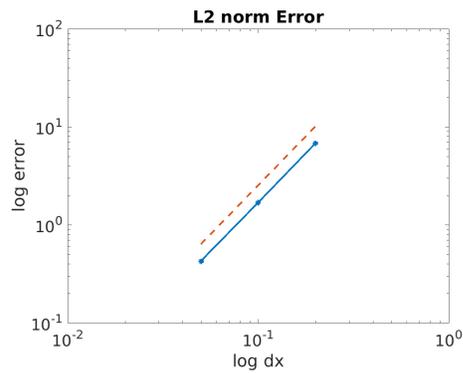

Fig. 9: Torus Scale Invariance

**6. Multidimensional Scaling Plots.** It is not clear how to visualize the similarities between surfaces if we are given the vectors of Shape DNA's. Therefore, we examine the similarities through multidimensional scaling. Multidimensional scaling is a non-supervised learning algorithm that takes in an object with many features and clusters them through nonlinear dimension reduction. It does this through minimizing a cost function:

$$f(Z) = \sum_{i=1}^{n} \sum_{j=1}^{n} (d_2(z_i - z_j) - d_1(x_i - x_j)). \tag{10}$$

This algorithm directly optimizes locations of objects in the $z$ space (in the space we want) from locations of objects in the original space $x$. We can use different metrics each step. $d_1$ is the high dimensional distance we want to match in the original space, $d_2$ is the low dimensional distance we can control, which is the plot distance and $d_3$ controls how we compare the objects in the two spaces. The result is the optimized objects locations for an N dimensional plot.

Once we compute the spectra of different surfaces, we can compute the dissimilarity matrix $D$ which contains pairwise distances between every spectrum. In other words, $D$ contains the distances between objects in the original space. For the results below, the standard Euclidean distances were used, which represents the metric Then, the MATLAB command 'mdscale' runs the multidimensional scaling algorithm and yields locations on the 2D or 3D plot such that the similarity distances between objects is preserved on the plot.

**6.1. Torus Experiment.** We test a very simple similarity like the thickness of the tori. We expect tori of similar thickness to be clustered together in the MDS plot. Torus 1 is the thickest torus, and torus 6 is the thinnest. Major radius was fixed to 1 and minor radii were 0.4, 0.3, 0.2, 0.1, 0.15, 0.1, 0.05. The shapes were resolved using grid size $\Delta x = 0.05$ such that the thinnest torus can be properly resolved by the Closest Point method.

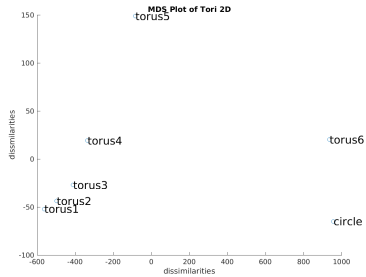

Fig. 10: Torus MDS Plot

Between 'torus 1', 'torus 2', 'torus 3', and 'torus 4', the minor radius decrements by 0.1 where as for the rest, it decrements by 0.05. However, we observe a more distinct clustering for the thicker tori. Our hypothesis is that the thickness of the tori for 5 and 6 is still poorly resolved by the given grid size, $\Delta x = 0.05$. Also note that the vertical dissimilarity distance scale is smaller than the horizontal scale by magnitudes. While it may seem that 4 and 5 are oddly far apart, their horizontal distance is reasonable. 6 is closest to the circle, which is an expected clustering. This experiment demonstrates that the MDS plots are able to capture thickness of objects.

**6.2. Holes on Sphere Experiment.** We test various sizes of holes on a unit sphere. Sphere ring 1 and 2 are spheres with small punctures with radii 0.05 and 0.1. Sphere ring 4 and 5 are big holes with radii 0.90 and 0.95, close to that of a unit hemisphere. Sphere ring 3 has a medium sized hole with radius 0.5. Ring is a sphere cut both top and bottom with radius 0.3 and 0.5 respectively. Neumann boundary conditions were imposed on the open surfaces. Our hypothesis was that 1 and 2 will cluster with the sphere, 4 and 5 with the hemisphere, 3 in the middle and the ring with the mobius strip.

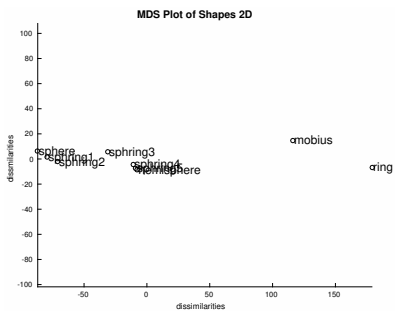

Fig. 11: Sphere Ring experiment

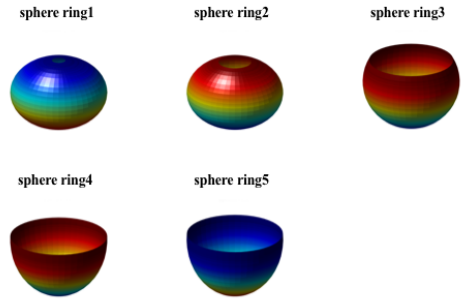

Fig. 12: Sphere Rings

In this plot, the MDS better recognizes global similarities than particularities. For example, the general shape of a sphere, hemisphere and ring is distinguished whereas the differences in the openness and the orientation of objects are more subtle.

**6.3. 2D vs. 3D.** We demonstrate the potential effectiveness of a 3D plot. With higher dimensional plots, the optimized location of each surface in the corresponding space contains more information about the Shape DNA's, although it would be harder to visualize. Our hypothesis is that adding a third axis to the 2D plot will capture additional characteristics of the surfaces as they can vary along more axes, with each axis representing a characteristic. The tested shapes are two sphere rings with Dirichlet and Neumann conditions, two hemispheres with Dirichlet and Neumann conditions, a sphere and an apple. These surfaces differ by their openness, boundary conditions and the general geometry.

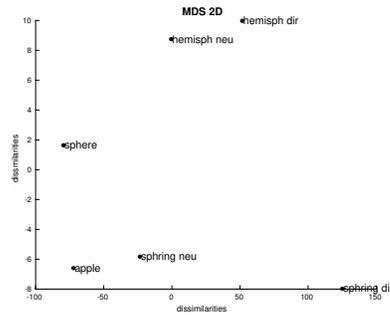

In the 2D plot, the horizontal axis scale is much more significant. We observe that the closed surfaces (apple and sphere) are at the left and the open surfaces scatter to the right. However, the boundary conditions are not distinguished and it is unclear how the surfaces vary along the vertical axis. The goal is achieved by the 3D plot as each axis shows a unique type of characteristic that the surfaces vary upon.

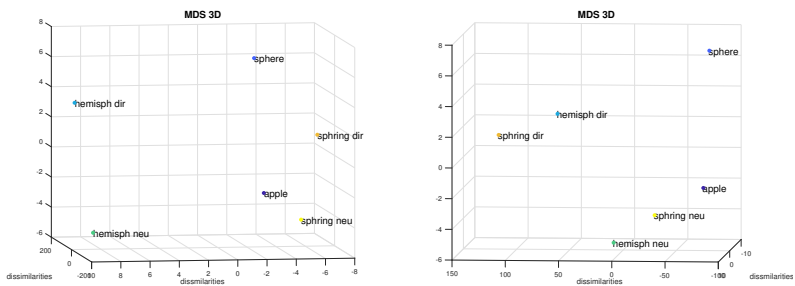

Fig. 13: The $x$-axis distinguishing the geometry

Fig. 14: The y-axis distinguishing the boundary conditions

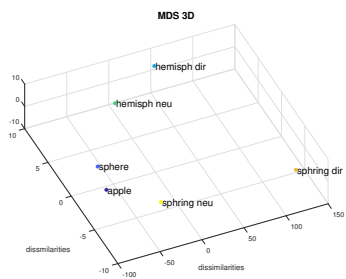

Fig. 15: The z-axis distinguishing the open and closed surfaces

**6.4. Random Shapes.** Finally, we compute the spectra of seemingly unrelated objects and try to uncover similarities using MDS. Note that the $x$-axis and $y$-axis are scaled differently, so the horizontal distance on the plot is greater than the vertical distance.

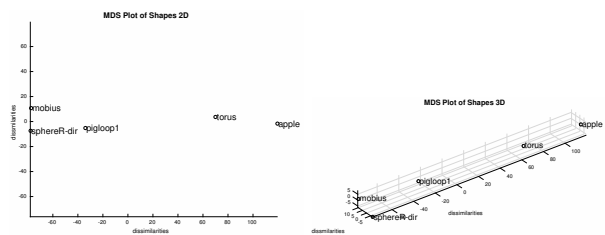

Fig. 16: MDS plots of Random Shapes

The hypothesis given this collection of shapes was that the objects with bound-aries are clustered (mobius strip and sphere ring), round objects (torus, apple) are clustered, and the triangulated pigloop2 (pig that has been smoothed by two Loop subdivisions) stands alone, in between the two clusters. To confirm this hypothesis, we add more triangulated objects (a non-smooth pig–annie, pigloop1 that has been smoothed by one Loop subdivision and a smooth bunny). We also add a sphere to verify that it clusters with the closed, round objects.

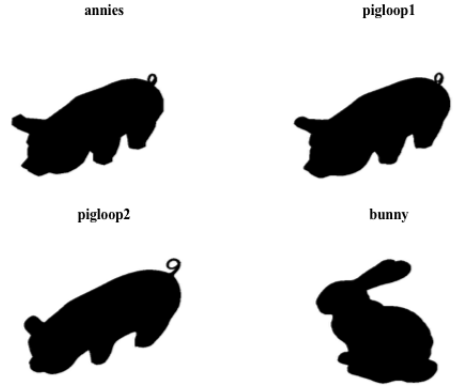

Fig. 17: Animals

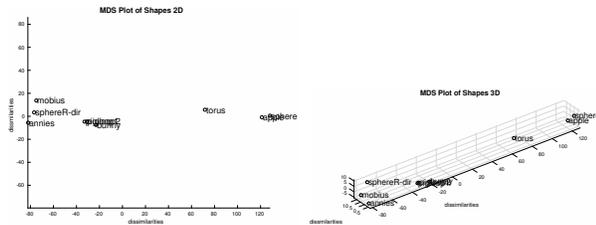

Fig. 18: MDS plots of Random Shapes

The result is as expected, except that the shape 'annies' is far from other animal shapes. One speculation is that since this particular pig has sharp edges and corners, it cannot be fully resolved by the Closest Point Method. Hence, it lies away from rest of the other pig shapes, and further away from the round objects. The other reasonably smooth triangulated animals lie in between. The sphere that was added clusters nicely with the round objects.

**Acknowledgments.** This work was done through Undergraduate Summer Research Award from NSERC, under the supervision of Dr. Colin Macdonald. I would like to thank Colin for his invaluable support and guidance on this project, as well as Dr. Steve Ruuth for motivating this work and teaching the Closest Point method during 2017 PIMS-SFU Summer School, and Chingyi Tsoi for her helpful contribution to the project. Also, many thanks to the UBC Math IT Department for the computational support.

## REFERENCES

[1] C. Gordon, D. Webb, and S. Wolpert, *Isospectral plane domains and surfaces via Riemannian orbifolds.*, Inventiones Mathematicae, 110 (1992), p. 1, https://doi.org/10.1007/BF01231320.
[2] M. Kac, *Can one hear the shape of a drum?*, American Mathematical Monthly, 73 (1996).
[3] R. J. LeVeque, *Finite Difference Methods for Ordinary and Partial Differential Equations*, SIAM, Philadelphia, PA, 2007.
[4] C. B. Macdonald, J. Brandman, and S. J. Ruuth, *Solving eigenvalue problems on curved surfaces using the closest point method*, J. Comput. Phys., 230 (2011).
[5] M. Reuter, S. Biasotti, G. P. D. Giorgi, and M. Spagnuol, *Discrete Laplace–Beltrami operators for shape analysis and segmentation*, Computers & Graphics, 33 (2009).
[6] M. Reuter, F.-E. Wolter, and N. Peinecke, *Laplace–Beltrami spectra as 'shape-DNA' of surfaces and solids*, Computer-Aided Design, 38 (2006).
[7] S. J. Ruuth and B. Merriman, *A simple embedding method for solving partial differential equations on surfaces*, J. Comput. Phys., 227 (2008).
[8] S. Seo, M. Chung, and H. Vorperian, *Heat kernel smoothing using Laplace–Beltrami eigenfunctions*, Medical Image Computing and Computer-Assisted Intervention, 230 (2010).
[9] I. von Glehn, T. März, and C. B. Macdonald, *An embedded method-of-lines approach to solving partial differential equations on surfaces*, ArXiv e-prints, (2013), https://arxiv.org/abs/1307.5657.